\documentclass{article}
\usepackage[utf8]{inputenc}
\usepackage{amsmath}
\usepackage{amssymb}
\usepackage{siunitx}
\usepackage{graphicx}
\usepackage{hyphenat}
\usepackage{authblk}
\usepackage{amsthm}
\usepackage[greek,english]{babel}
\usepackage{lmodern}
\usepackage[numbers]{natbib}

\title{A note on the control \\ 
of processes exhibiting input multiplicity}
\author[1]{Robert J. Lovelett}
\author[2]{Yorgos M. Psarellis}
\author[1,2]{Ioannis G. Kevrekidis}
\author[3,*]{Manfred Morari}
\affil[1]{Department of Chemical and Biological Engineering, Princeton University}
\affil[2]{Department of Chemical and Biomolecular Engineering, Johns Hopkins University}
\affil[3]{Department of Electrical and Systems Engineering, University of Pennsylvania}
\affil[*]{Correspondence: morari@seas.upenn.edu}
\begin{document}
\maketitle

\theoremstyle{definition}
\newtheorem*{definition}{Definition}
\newtheorem*{theorem}{Theorem}



\begin{abstract}
Steady state multiplicity can occur in nonlinear systems, and this presents challenges to feedback control. {\em Input multiplicity} arises when the same steady state output values can be reached with system inputs at different values. Dynamic systems with input multiplicities equipped with controllers with integral action have multiple stationary points, which may be locally stable or not. This is undesirable for operation. For a 2x2 example system with three stationary points we demonstrate how to design a set of two single loop controllers such that only one of the stationary points is locally stable, thus effectively eliminating the ``input multiplicity problem'' for control. We also show that when MPC is used for the example system, all three closed-loop stationary points are stable. Depending on the initial value of the input variables, the closed loop system under MPC may converge to different steady state input instances (but the same output steady state). Therefore we computationally explore the basin boundaries of this closed loop system. It is not clear how MPC or other modern nonlinear controllers could be designed so that only specific equilibrium points are stable.
\end{abstract}

Dedication:
RJL completed his PhD under Tunde Ogunnaike's supervision at the University of Delaware.
He is grateful that despite Tunde's many roles---professor, dean, author, journal editor, and more---Tunde reserved time and attention to support his students' lives and careers.
He will miss Tunde's kind mentorship and thoughtful guidance, and the opportunity to discuss and share new research findings and life updates.
IGK (who will miss Tunde's trademark Greek greeting ``\textgreek{παλιοφιλε},'' old friend) kept in contact with him from the DuPont days to University of Delaware.
His last invitation to Tunde for a seminar at Hopkins was initially delayed by Covid, never to materialize.
MM was on the faculty at the University of Wisconsin when Tunde was there for his graduate studies. 
In his typical kind and inquisitive manner Tunde contributed to MM's research about the Relative Gain Array (RGA) on which the ideas in this note are based
It was the starting point for many stimulating exchanges in the following 40 years. 

\section{Introduction}

In this note, we are studying the control of nonlinear time invariant square systems with manipulated input $\mathbf{u} \in \mathbb{R}^n$ and controlled output $\mathbf{y} \in \mathbb{R}^n$. The system is assumed to be stable in the sense that any asymptotically constant input $\mathbf{u}$ will result in an asymptotically constant output $\mathbf{y}$. Thus the system defines a static map $\mathbf{G}$ relating the constant input vector $\mathbf{u}$ to the eventually constant output vector $\mathbf{y}$:
\begin{equation}
\mathbf{y}=\mathbf{G}(\mathbf{u}).
\end{equation}

\begin{figure}
    \centering
    \includegraphics[width=.5\textwidth]{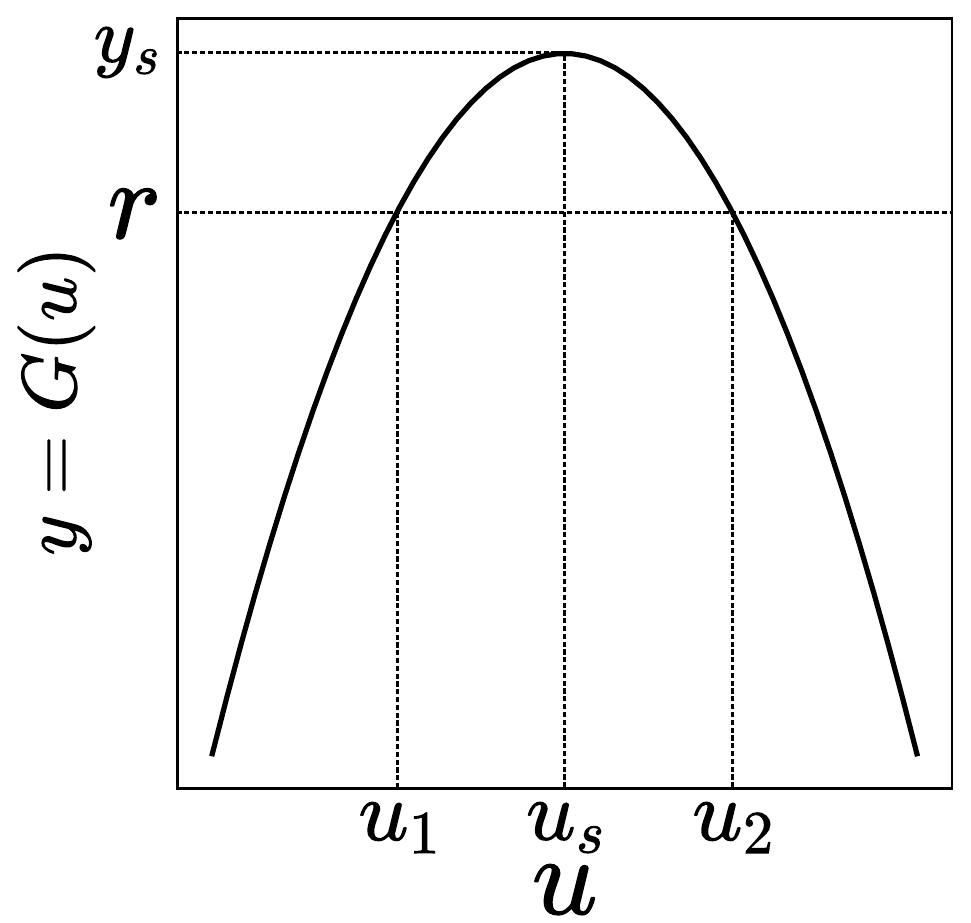}
    \caption{A simple SISO system that exhibits input multiplicity. For setpoint $r$, there are two values of $u$ ($u_1$ and $u_2$) for which $y=G(u)=r$. For smooth $G$, a necessary condition for this phenomenon is that $G'(u)=0$ for some value of $u$.}
    \label{fig:SISO_example}
\end{figure}

We are interested in systems where the gain matrix, i.e., the Jacobian $\mathbf{G}'(\mathbf{u})= \left. \frac{ \partial\mathbf{G} }{ \partial\mathbf{u} }\right|_{\mathbf{u}}$ is singular for some $\mathbf{u}_s$, i.e. $\mathrm{det} \left. \left( \frac{ \partial \mathbf{G}} {\partial \mathbf{u}}\right|_{\mathbf{u_s}} \right) =0$.
A simple SISO  example is shown in Fig. \ref{fig:SISO_example}.
The map $\mathbf{G}(\mathbf{u})$ has either an extremum or a saddle point at $\mathbf{u}=\mathbf{u}_s$.
In either case it is impossible to control the system at a reference $\mathbf{r}=\mathbf{y}_s$ because the controller gain vanishes.
In the MIMO case, the condition $\mathrm{det} \left. \left( \frac{ \partial \mathbf{G}} {\partial \mathbf{u}}\right|_{\mathbf{u}} \right) =0$ implies that the gain vanishes in a certain input direction, again indicating that control is not possible. 

The output variable $\mathbf{y}$ may be associated with process economics so that its maximization is sought.
In this case extremum seeking control schemes would be used, which have been studied extensively \cite{ariyur2003real} and are not the subject of this note.
Here we study only regulatory control, i.e. control of the system at a specified setpoint $\mathbf{y}=\mathbf{r}$.

It should be emphasized that this lack of (practical) controllability at $\mathbf{r}=\mathbf{y}_s$ is a property of the system and cannot be fixed by any clever controller design.
Instead, it is necessary to change the system itself, for example, by choosing a different manipulated variable or by changing the design parameters of the system such that the singular point is pushed out of the region of interest.

The system may be difficult to control even at setpoints $\mathbf{r} \neq \mathbf{y}_s$ where the gain does not vanish. 
In Fig. \ref{fig:SISO_example} we observe ``input multiplicity'', i.e. that a particular desired output $y=G(u)=r$ can be reached for multiple inputs $u$.
In practical operation it may not be easily predictable which input will be ultimately ``chosen'' by the controller.
This could be undesirable because it would make it very difficult to monitor the process for ``correct'' operation.

It appears that the control issues surrounding input multiplicity were first discussed by Koppel \cite{Koppel1982}.
In this pioneering paper and several follow-ups \cite{Koppel1983, Dash1989} he presents various examples, which suggest that input multiplicity can occur commonly in process control.
The fact that it has not been discussed in the literature more broadly in the past several decades may simply mean that the erratic behavior associated with input multiplicity has been observed in practical operation, but not diagnosed correctly.
Such a focused analysis would have required an accurate process model which may not have been available.

In our note we will analyze two different approaches to design controllers for systems with input multiplicities. 

\begin{enumerate}
\item{Multi-loop integral control.  Our control objective is slow control which eliminates off-set. We will achieve this by a set of pure integral controllers, which guarantee local asymptotic stability for a specific input instance. }

\item{Model Predictive Control (MPC). We determine the control action by minimizing a quadratic objective over a moving finite horizon. If the closed loop system is locally stable, then there will be no off-set for our particular choice of objective.}
\end{enumerate}

\section{Multi-loop integral control}

Based on his experience, Koppel argues convincingly that in practical process control the open loop speed of response is usually satisfactory, we rarely wish to make it faster.
Even if we wanted to, any speed-up would require a more accurate dynamic process model, which is typically not available.
We do want off-set free tracking, however, which we can simply achieve with integral control action.  

In this part of our note we follow this premise.
Our goal is not to demonstrate that some elaborate nonlinear controller assures good dynamic performance for a specific example process.
We want to derive general theoretical results which are useful for the analysis and design of controllers for open-loop stable processes for which the gain matrix becomes singular within the domain of operation and where the performance specifications require integral control action. 

Just like Koppel (and Bristol \cite{Bristol1966} and Niederlinski \cite{Niederlinksi1971} whom he references) our theoretical results require only knowledge of the static model of the system captured in the nonlinear map $\mathbf{G}$.
Somewhat surprisingly, nobody (us included) seems to have connected the work by Koppel to the works by Grosdidier et al. \cite{Grosdidier1985} and Campo and Morari \cite{Campo1994} which are addressing this problem in exactly this manner in a precise theoretical framework. 

Let us assume that for an output reference $\mathbf{y}=\mathbf{r}$, the static map $\mathbf{G}$ allows several inputs $\mathbf{u}_1$, $\mathbf{u}_2$, \ldots,  i.e., 
\begin{equation}
\mathbf{r}=\mathbf{y}=\mathbf{G}(\mathbf{u}_1)=\mathbf{G}(\mathbf{u}_2)=\ldots
\end{equation}
Then we will derive necessary conditions to be met by an input $\mathbf{u}_i$ such that a dynamic controller with integral action exists which assures local closed loop properties like stability and fault tolerance. 

Consider, for example, the SISO case in Fig. \ref{fig:SISO_example}, where the two inputs $u_1$ and $u_2$ both yield the output $y=r$.
Our results establish that for steady state $(y,u_1)$ there exists a controller with integral action which stabilizes the system at that state. 
With this controller the other steady state $(y,u_2)$ is unstable.
The reverse is true for any controller stabilizing $(y,u_2)$.
This is important---the second potential equilibrium point arising from input multiplicity is \emph{always} unstable and therefore effectively irrelevant for controller design and system operation.
This SISO example is obvious, just reinterpreting the well-known fact that in the presence of integral control, stability requires the sign of the process gain not to change.
The other parts of the controller and the system dynamics do not matter.

In the MIMO case the issues are more complex.
In the rest of this section we will give some indication how they could be addressed.
Again, as in the SISO case, our objective for practical operation will be to assure that through appropriate controller design, there is only a \emph{single} stable state.
For this purpose, we first need to derive a MIMO generalization of the SISO positive gain condition we mentioned above.
This was done in great detail in the paper by Grosdidier et al. \cite{Grosdidier1985} for linear systems.
Even  more details have been reported by Campo and Morari \cite{Campo1994}. Here we will recall the main definitions and results. Later we will use these results to analyze the local stability of the nonlinear systems studied in this note.

\paragraph{Integral Controllers}

\begin{figure}
    \centering
    \includegraphics[width=.75\textwidth]{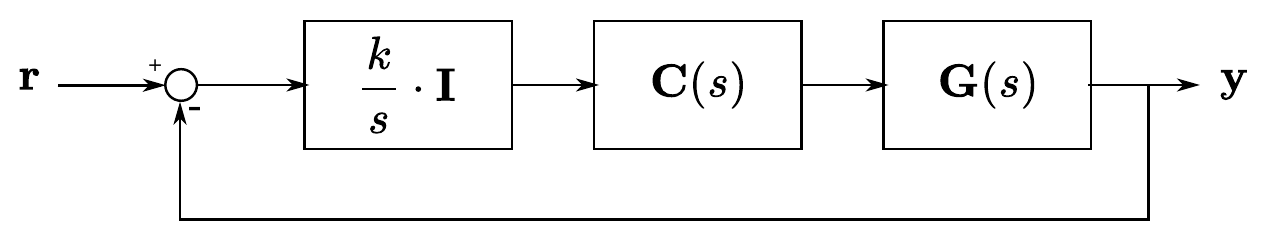}
    \caption{Basic MIMO integral control configuration.}
    \label{fig:mimo_config}
\end{figure}

In the present context, the term integral controller is used to designate PI, PID controllers, or any multi-variable feedback controller which includes integral action.
All such controllers can be decomposed into a matrix of integrators $\frac{k}{s}\cdot\mathbf{I}$ and a compensator matrix $\mathbf{C}(s)$.
Here, $k$ is a positive constant and $\mathbf{I}$ is the identity matrix.
Such a decomposition simplifies the analysis since we are not interested in a specific compensator but rather in the general consequences of the integral action.
Next, we consider the $n \times n$ plant transfer matrix $\mathbf{G}(s)$ and the control configuration shown in Figure \ref{fig:mimo_config}.
Define $\mathbf{H}(s) = \mathbf{G}(s)\mathbf{C}(s)$.

\begin{definition}[Integral Controllability]
The open-loop stable system $\mathbf{H}(s)$ is called integral controllable if there exists a $k^* > 0$ such that the closed-loop system shown in Fig. \ref{fig:mimo_config} is stable for all values of $k$ satisfying $0 < k < k^*$ and has zero tracking error for asymptotically constant disturbances.
\end{definition}

\begin{theorem} The rational system $\mathbf{H}(s)$ is integral controllable if all the eigenvalues of $\mathbf{H}(0)$ lie in the open right half complex plane.
The rational system $\mathbf{H}(s)$ is not integral controllable if any of the eigenvalues of $\mathbf{H}(0)$ lie in the open left half complex plane.
(The Theorem says nothing about systems for which the eigenvalues of $\mathbf{H}(0)$ lie in the closed right half plane and include eigenvalues on the imaginary axis.)
\end{theorem}

Note that this theorem provides the complete generalization of the positive gain condition for SISO systems stated earlier. In his paper Koppel requires $\mathrm{det}\left(\mathbf{G}(0)\right)>0$ which follows directly.

We assume that $\mathbf{C}(s)$ is the product of a diagonal matrix with a permutation matrix, i. e., we assume that we can decide on the variable pairing (permutation matrix) and then choose the remaining part of the controller to be diagonal, i.e., a set of single loop controllers. 

For the MIMO control of systems with input multiplicity we want to choose the controller $\mathbf{C}$ such that for only one of the possible stationary points all the eigenvalues of the linearized system $\mathbf{H}$ are in the RHP.
This would imply that the closed loop system has a single stationary stable point which would be attractive for operation.
There may be many controllers $\mathbf{C}$ which lead to the desired sets of eigenvalues.
A numerical search for those that also lead to good control performance may be difficult.
Therefore, we  will construct the multi-loop SISO controllers via sequential loop closing. 
We will first decide on the pairings, then on the sequence of closing the loops. 
In this manner, for each loop we can easily determine the sign of the loop gain required for stability from the product of the system gain and the appropriate element of the RGA.

\subsection{MIMO Case Study}

\begin{figure}
    \centering
    \includegraphics[width=.75\textwidth]{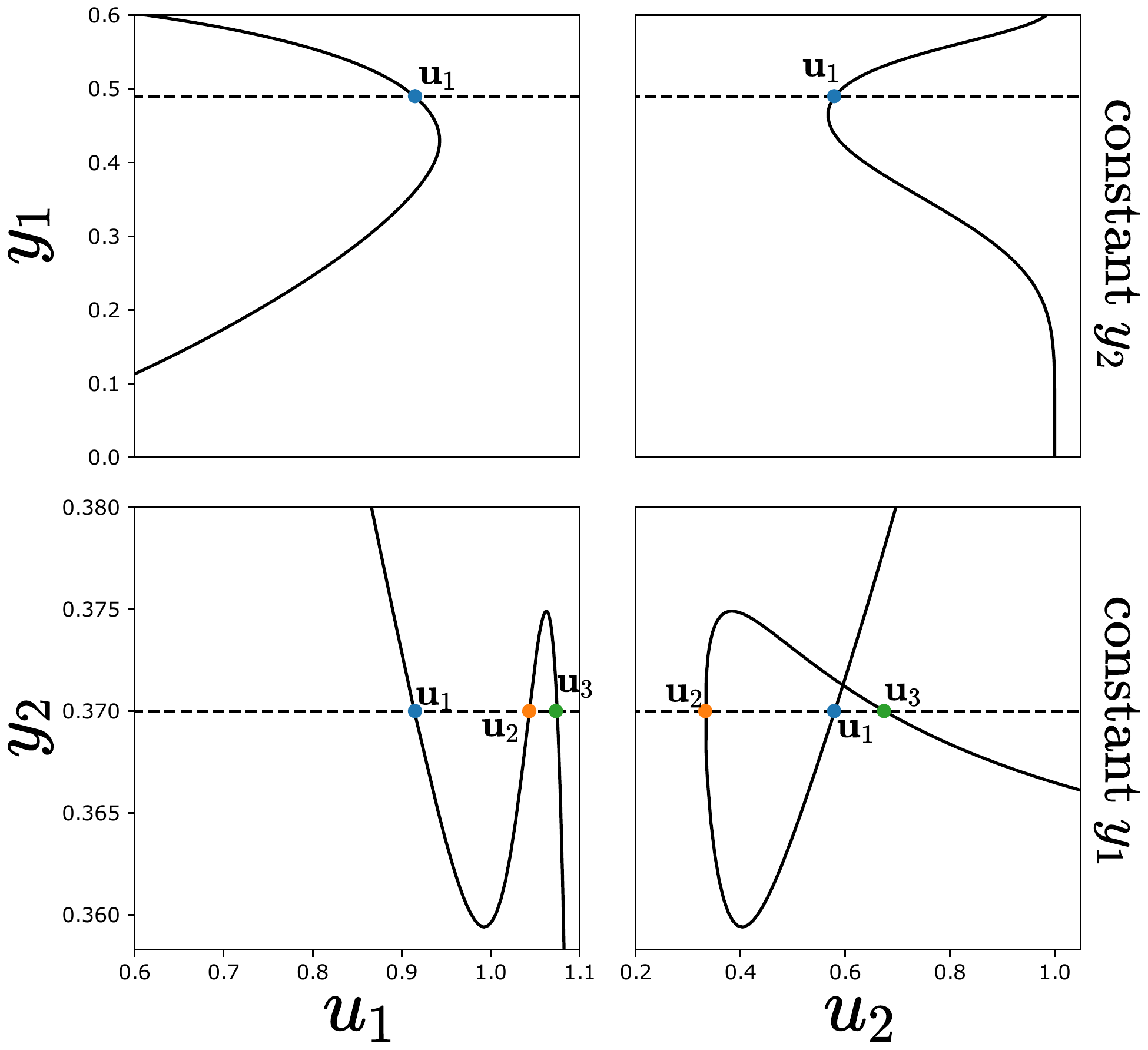}
    \caption{Numerical continuation of \emph{inputs} ($\mathbf{u}$) with respect to the \emph{outputs} ($\mathbf{y}$). Columns correspond to the two inputs, while rows to the output for which the continuation was performed (with the other one kept fixed). The system was initialized at input instance 1, $\mathbf{u}_1$. We caution the reader that the sign of the open loop gain $\mathbf{G}(s)$ \emph{cannot} simply be observed from slopes at each steady state on this figure because the \emph{other} output variable is held constant for numerical continuation; it is the \emph{closed loop} gain that is observed.}
    \label{fig:cont}
\end{figure}

We will demonstrate how to construct the multi-loop SISO controllers by way of a case study of a $2\times2$ system with three stationary points.
We use the example from \cite{Koppel1982} detailed in the Appendix, in which two sequential reactions take place in a continuous stirred tank reactor, where inputs $\mathbf{u}$ are scaled temperature and residence time and outputs $\mathbf{y}$ are mole fractions of two (out of three) chemical species.
Figure \ref{fig:cont} shows results from numerical continuation of each of the \emph{inputs}, starting from steady state $\mathbf{u}_1$; each of the three open-loop stable steady states, $\mathbf{u}_1, \mathbf{u}_2, \mathbf{u}_3$, are observed.
The results of our multi-loop SISO controller construction are summarized in Table \ref{tab:icontrol}, where we determined the signs of the gains of the diagonal elements of $\mathbf{C}$ via sequential loop closing starting with Pairing 1 following Sequence 1 etc..

In conclusion for each equilibrium point there is at least one pairing and gain sign choice such that the selected equilibrium is the only local stable one. 
\begin{itemize}
\item{
For Pairing $(y_1-u_1), (y_2-u_2)$
with gains  ++, only $\mathbf{u}_2$ is locally stable.
With gains +\textminus~only, $\mathbf{u}_3$ is locally stable.}

\item{
For Pairing  $(y_1-u_2), (y_2-u_1)$
with gains  ++, only $\mathbf{u}_1$ is locally stable.
With gains  \textminus+, only $\mathbf{u}_2$ is locally stable.
}
\end{itemize}

This demonstrated simple controller design to achieve a single asymptotically stable equilibrium point may be attractive for practical operation.
Note, however, that other dynamic phenomena like limit cycles could be present. 
It would be interesting to explore whether this is possible for the simple system studied here, comprising a static process with single loop integral controllers and if it can be avoided by choosing the integral gain sufficiently small as assumed for the definition of integral controllability. 

\begin{table}
\centering
\caption{Gain matrix, relative gain array, and signs of integral gains that stabilize each input setting of the system from Equation \ref{eqn:system}, for each controller configuration (signs are ordered by the \emph{controlled} variables, e.g., ``\textminus+'' indicates that the system may be stable when the $y_1$ controller has negative gain and the $y_2$ controller has positive gain).}
\label{tab:icontrol}
\makebox[\linewidth]{
\begin{tabular}{cccccc}
\hline
$\mathbf{y}=\mathbf{r}$ &
  $\mathbf{u}$ &
  $\mathbf{G}$ &
  RGA &
  \begin{tabular}[c]{@{}l@{}}$(y_1-u_1),$ \\ $(y_2-u_2)$\end{tabular} &
  \begin{tabular}[c]{@{}l@{}}$(y_1-u_2),$\\ $(y_2-u_1)$\end{tabular} \\ \hline
$(0.49, 0.37)$ &
  $(0.914, 0.580)$ &
  \begin{tabular}[c]{@{}l@{}}$\begin{bmatrix}\\ -1.89 & -0.78 \\\\  0.17 & 0.15\\ \end{bmatrix}$\end{tabular} &
  \begin{tabular}[c]{@{}l@{}}$\begin{bmatrix}\\  1.85& -0.85 \\\\ -0.85&  1.85\\ \end{bmatrix}$\end{tabular} &
  \textminus+ &
  \textminus\textminus,++ \\
$(0.49, 0.37)$ &
  $(1.043, 0.333)$ &
  \begin{tabular}[c]{@{}l@{}}$\begin{bmatrix}\\ 0.037& -0.805 \\\\  0.361 &  0.256\\ \end{bmatrix}$\end{tabular} &
  \begin{tabular}[c]{@{}l@{}}$\begin{bmatrix}\\ 0.032& 0.968 \\\\ 0.968& 0.032\\ \end{bmatrix}$\end{tabular} &
  ++ &
  \textminus+ \\
$(0.49, 0.37)$ &
  $(1.07, 0.674)$ &
  \begin{tabular}[c]{@{}l@{}}$\begin{bmatrix}\\ 14.0& -0.217 \\\\ -7.24 &  0.097\\ \end{bmatrix}$\end{tabular} &
  \begin{tabular}[c]{@{}l@{}}$\begin{bmatrix}\\ -6.58& 7.58 \\\\ 7.58& -6.58\\ \end{bmatrix}$\end{tabular} &
  +\textminus,\textminus+ &
  \textminus\textminus \\ \hline
\end{tabular}
}
\end{table}

\section{Model Predictive Control}

We continue our case study of the MIMO CSTR system by examining a model predictive control scheme.
Unlike the case of integral controllers, where we can guarantee by construction that only one of the three steady state input settings is stable, the model predictive controller does not provide any such assurances.
Consider the nonlinear MPC controller, defined by the solution to the optimization problem
\begin{equation}\label{eqn:mpc_problem}
\begin{aligned}
\min_{\mathbf{u}} F(\mathbf{u}) = \sum_{k=0}^H
\left(
    \left((\mathbf{y}(k+1)-\mathbf{r}) \mathbf{K}_y (\mathbf{y}(k+1)-\mathbf{r})^T \right) 
    +
    \left((\mathbf{u}(k+1)-\mathbf{u}(k)) \mathbf{K}_u (\mathbf{u}(k+1)-\mathbf{u}(k))^T \right)
    \right)\\
\end{aligned}
\end{equation}
where $\mathbf{y}(k)$ and $\mathbf{u}(k)$ are sampled time points from the system with sampling time of 0.5, and the tuning parameters are:
\begin{equation*}
    \begin{aligned}
    H&=2\\
    \mathbf{K}_y &= \mathbf{I}\\
    \mathbf{K}_u &= 2\mathbf{I}
    \end{aligned}
\end{equation*}

Note that in this formulation, there is a penalty (the term with $\mathbf{K}_u$) on the change in input variables; this objective function assures that there will be off-set free regulatory control if the system is locally stable.
In general, the optimization problem is not convex and requires a global solver.
In practice, a local solution is obtained using the BFGS algorithm (as implemented in SciPy \cite{2020SciPy-NMeth}) and with initial guess $\mathbf{u}=\begin{bmatrix} \mathbf{u}_0 & \mathbf{u}_0 \end{bmatrix}$, where $\mathbf{u}_0$ is the vector of inputs at the previous time step.

\begin{figure}
    \centering
    \includegraphics[width=0.9\textwidth]{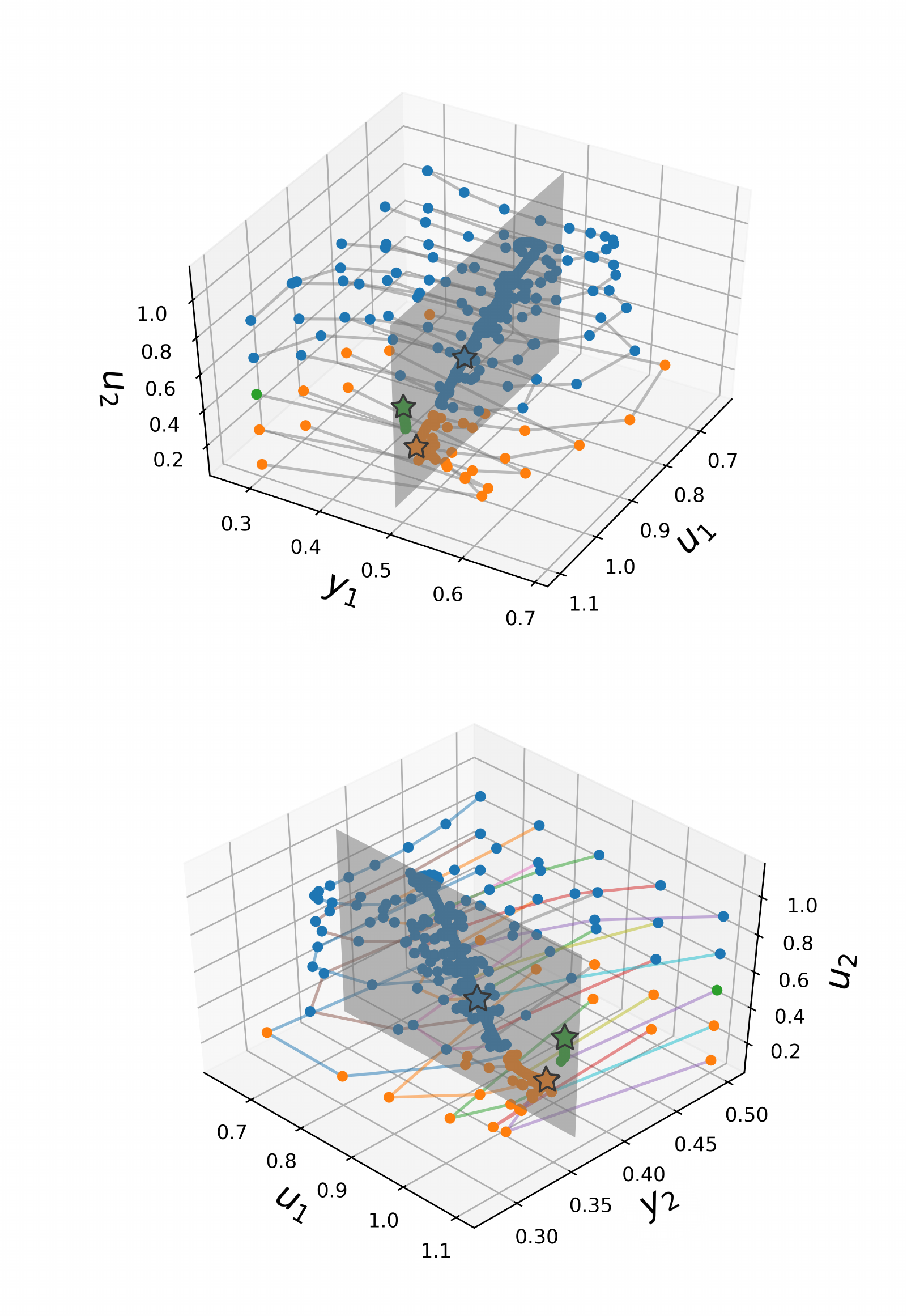}
    \caption{Trajectories of the closed loop MPC system, showing 3D projections of the 4D system in $(u_1, y_1, u_2)$ (top) and $(u_1, y_2, u_2)$ (bottom). The set point and input instances, which comprise the fixed points, are denoted by the ``star'' points. Trajectories are colored by the fixed points to which they are attracted. The shaded plane shows the set points for $y_1$ (top) and $y_2$ (bottom), indicating that all fixed points lie at the set point.}
    \label{fig:MPC_trajectories}
\end{figure}

Figure \ref{fig:MPC_trajectories} shows 3D projections of trajectories of the 4D MIMO system using the MPC controller at various initial input ($\mathbf{u}$) conditions, given the same initial position in the controlled variables ($\mathbf{y}_0=(0.3,0.5)$).
We observe that for each initial condition, the system converges to the set point, but at different values of the input variables.
In all three cases, the system exhibits time scale separation, where the system approaches a slow manifold in several time steps, and then very slowly approaches the steady state set point.
These results suggest that the system has three basins of attraction that demarcate which fixed point the closed loop system approaches at steady state.
Each basin of attraction contains a single fixed point, whose position in the closed loop system each corresponds to one of the three input instances in Equation \ref{eqn:steady_states} in the Appendix. 

\begin{figure}
    \centering
    \includegraphics[width=\textwidth ]{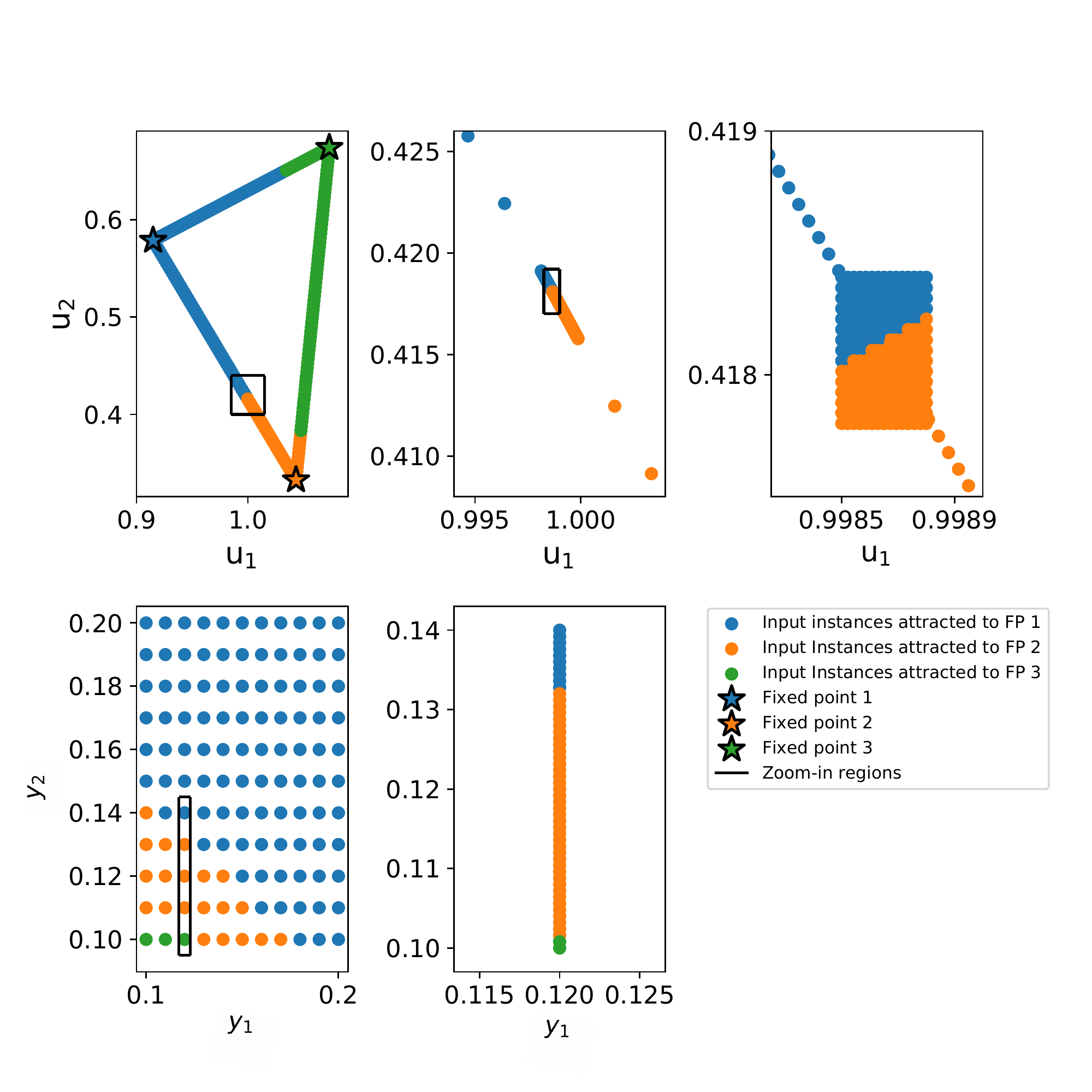}
    \caption{Each point is an initial condition in $\mathbf{u}$-space (top row) or $\mathbf{y}$-space (bottom row) and is colored by the fixed point in which the closed loop system will converge towards. For all points on the top row, $\mathbf{y}_0=(0.3, 0.5)$; and for all points on the bottom row, $\mathbf{u}_0=(0.96, 0.5)$}.
    \label{fig:zoomed_basins}
\end{figure}

We attempted to characterize the basins of attraction of the closed loop system under MPC.
In nonlinear systems, the boundaries of basins of attraction may be smooth or they may be fractal \cite{McDonald1985, grebogi1987chaos, Thompson_1988}.
These basin boundaries are often associated with saddle-type invariant objects (e.g. stable manifolds of saddle steady states/fixed points). 
We used numerical experiments to examine whether the basin boundaries in our system are smooth or fractal, and visualized the results with 2D projections of the 4D phase space. 
Figure \ref{fig:zoomed_basins} shows initial conditions in $\mathbf{u}$-space (top row) and $\mathbf{y}$-space (bottom row) colored by the steady state to which they are attracted; in particular, we visualize regions very close to the basin boundaries (in the far right plots).
In both the $(u_1,u_2)$ plane and $(y_1,y_2)$ plane, we do not observe any evidence of a fractal basin of attraction, suggesting that the basin boundaries are, to the best of our ability to discern, smooth.
We were not able to numerically locate any additional saddle-type, unstable invariant objects that organize these basin boundaries. 

\section{Discussion}

Input multiplicity presents a challenge for control system design.
Here, we compare controller design by sequential single loop closing using integral controllers and a nonlinear model predictive control approach.
By using the simple integral controller, we demonstrate that we can guarantee by construction that the closed loop system has \emph{only one} locally stable fixed point; alternatively the closed loop system under model predictive control has three locally stable fixed points that make it challenging to predict which fixed point the system will approach.
The basin boundaries that separate the fixed points in the MPC case, although smooth in our model system, may be especially challenging to identify in practice, especially in the face of plant-model mismatch.
These considerations favor the use of the integral control system, where we can assure local stability and off-set free tracking by construction.

\bibliographystyle{abbrvnat}
\bibliography{refs}

\appendix
\section{MIMO CSTR System Description}
Our case study is a bivariate system that exhibits input multiplicity, which was first examined by Koppel \cite{Koppel1982}.
The system describes a CSTR with two sequential chemical reactions and three chemical species ($x_1$, $x_2$, and $x_3=1-x_1-x_2$; $x_3$ is redundant and is not used in the model).
The state equations are:
\begin{equation}\label{eqn:system}
\begin{aligned}
    f_1(\mathbf{x}) &= -k_1 x_1 + k_4 x_2 + \frac{(x_{10} - x_1)}{\tau}\\
    f_2(\mathbf{x}) &=  k_1 x_1 - k_2 x_2 + k_3 (1-x_1-x_2) + \frac{(x_{20}-x_2)}{\tau},
\end{aligned}
\end{equation}
where: 
\begin{equation}
\begin{aligned}
k_1 &= k_{10}\exp\left(-\frac{E_1}{RT_0} (T_0/T - 1)\right)\\
k_2 &= k_{20}\exp\left(-\frac{E_2}{RT_0} (T_0/T - 1)\right)\\
k_3 &= k_{30}\exp\left(-\frac{E_3}{RT_0} (T_0/T - 1)\right)\\
k_4 &= k_{40}\exp\left(-\frac{E_4}{RT_0} (T_0/T - 1)\right).
\end{aligned}
\end{equation}
The two inputs are scaled temperature and scaled residence time, defined as:
\begin{equation}
\begin{aligned}
u_1&= \frac{T}{T_0}\\
u_2&= \frac{\tau}{\tau_0+\tau}.
\end{aligned}
\end{equation}
The controlled variables are the non-redundant chemical species $x_1$ and $x_2$:
\begin{equation}
\mathbf{y} = \mathbf{h}(\mathbf{x}) = \mathbf{x}.
\end{equation}
The remaining, non-input, constant value parameters used by Koppel \cite{Koppel1982}, are given in Table \ref{tab:param}.
\begin{table}
\centering
\caption{System parameters.}
\label{tab:param}
\makebox[\linewidth]{
\begin{tabular}{ccc}
\hline
    Parameter        & Value \\
\hline
    $k_{10}$        & \SI{1.0}{\per \second} \\
    $k_{20}$        & \SI{0.7}{\per \second} \\
    $k_{30}$        & \SI{0.1}{\per \second} \\
    $k_{40}$        & \SI{0.007}{\per\second} \\
    $\frac{E_1}{R}$ & \SI{5000}{\kelvin} \\
    $\frac{E_2}{R}$ & \SI{6000}{\kelvin}\\
    $\frac{E_3}{R}$ & \SI{30000}{\kelvin}\\
    $\frac{E_4}{R}$ & \SI{50000}{\kelvin}\\
    $T_0$           & \SI{600}{\kelvin}\\
    $\tau_0$        & \SI{1}{\per \second}\\
    $x_{10}$        & 0.8\\
    $x_{20}$        & 0.2
    
\end{tabular}
}
\end{table}

Consider a multiple-input multiple output (MIMO) case, with $\mathbf{y}=\mathbf{h}(\mathbf{x})=\mathbf{x}$, and a set point that exhibits input multiplicity: $\mathbf{r} = (0.49, 0.37)$, with the following three possible steady state input instances:
\begin{equation}\label{eqn:steady_states}
    \begin{aligned}
    \mathbf{u}_{1}&=(0.914,0.580)\\
    \mathbf{u}_{2}&=(1.043,0.333)\\
    \mathbf{u}_{3}&=(1.07, 0.674)
    \end{aligned}
\end{equation}

All of the steady states are open-loop stable, having eigenvalues of the Jacobian matrix ($\mathbf{J}=\frac{\partial \mathbf{f}}{\partial \mathbf{x}}$) in the left half plane:
\begin{equation}
    \begin{aligned}
    \mathbf{\lambda}_{1} &= [\begin{matrix}-1.185 & -1.002\end{matrix}]\\
    \mathbf{\lambda}_{2} &= [\begin{matrix}-3.230 & -4.110\end{matrix}]\\
    \mathbf{\lambda}_{3} &= [\begin{matrix}-2.912 & -6.260\end{matrix}]
    \end{aligned}
\end{equation}

\end{document}